\documentclass[12pt]{article}
\newcommand{\di}{\displaystyle}
\usepackage{amsmath}
\usepackage{amssymb}
\usepackage{amsthm}
\newtheorem{proposition}{Proposition}
\newtheorem{remark}{Remark}
\newtheorem{definition}{Definition}
\newtheorem{theorem}{Theorem}

\begin{document}

\title{Rarita Schwinger Type Operators on Cylinders}

\author{Junxia Li and John Ryan\\
\emph{\small Department of Mathematics, University of Arkansas, Fayetteville, AR 72701, USA.}\\
Carmen J. Vanegas\\
\emph{\small Departamento de Matem\'aticas, Universidad Sim\'on Bol\'ivar, Caracas, Venezuela.}}

\date{}

\maketitle
\begin{abstract}

Here we  define Rarita-Schwinger operators on cylinders
and construct their fundamental solutions. Further the fundamental solutions to the cylindrical
Rarita-Schwinger type operators are achieved by applying translation groups.
In turn, a Borel-Pompeiu Formula, Cauchy Integral Formula and a Cauchy Transform are presented for the cylinders.
Moreover we show a construction of a number of conformally inequivalent spinor bundles on these cylinders.
Again we construct Rarita-Schwinger operators and their fundamental solutions in this setting. Finally we study the remaining Rarita-Schwinger type operators on cylinders.
\end{abstract}

\begin{center}
\small \bf {This paper is dedicated to the memory of Jaime Keller}
\end{center}

\section{Introduction}

The Rarita-Schwinger operators are generalizations of the Dirac operator which in turn is a natural
generalization of the Cauchy-Riemann operator. They have been studied in Euclidean space in \cite{BSSV, BSSV1,DLRV, Va, Va1} and on spheres and
real projective spaces in \cite{LRV}.

Conformally flat manifolds are manifolds with atlases whose transition functions are M\"{o}bius transformations. They can be
constructed by factoring out a subdomain $U$ of either the sphere $\mathbb{S}^n$ or $\mathbb{R}^n$ by a Kleinian subgroup $\Gamma$
of the M\"{o}bius group where $\Gamma$ acts strongly discontinuously on $U$. This gives rise to the conformally flat manifold
$U\setminus \Gamma$.
Real projective spaces are examples of conformally flat manifolds of type $\mathbb{S}^n\setminus \{\pm 1\}$.
Other simple examples are $\mathbb{R}^n\setminus \mathbb{Z}^l$, where $\mathbb{Z}^l$ is an integer lattice and
$1\leq l \leq n$. We call these manifolds cylinders.
In case $l=n$ we have the $n$-torus. For cylinders and tori the conformal structure is given by translations.

In this paper we  define Rarita-Schwinger operators on cylinders
and construct their fundamental solutions. Further the fundamental solutions to the cylindrical
Rarita-Schwinger type operators are achieved by applying translation groups.
In turn, Borel-Pompeiu Formula, Cauchy Integral Formula and a Cauchy Transform are presented for the cylinders.
We also show a construction of a number of conformally inequivalent spinor bundles on these cylinders.
Again we construct Rarita-Schwinger operators and their fundamental solutions in this setting.
The Rarita Schwinger type operators on tori are left for future research.

In \cite{KR} R. Krausshar and J. Ryan introduce Clifford analysis on cylinders and tori making use of the fact that the universal
covering space of all of these manifolds is $\mathbb{R}^n$. So provided the functions and kernels are $l$-periodic for
some $l \in \{1, \dots, n \}$ then R. Krausshar and J. Ryan use the projection map to obtain the equivalent function or kernel on those manifolds. Following them, we take the Rarita-Schwinger kernel in $\mathbb{R}^n$ and use the
translation group to construct new kernels that are $l$-fold periodic, which are projected on the cylinders.
In order to prove that the new kernels are well defined we adapt the Eisenstein series argument developed in \cite{K}. We also show a construction of a number of conformally inequivalent spinor bundles on these cylinders. In the last section, we discuss the remaining Rarita-Schwinger operators studied in \cite{LR} and determine their kernels on cylinders. We finally establish some basic integral formulas associated with the remaining Rarita-Schwinger operators on cylinders.

\section{Preliminaries}

\par A Clifford algebra, $Cl_{n},$ can be generated from $\mathbb{R}^n$ by considering the
relationship $$\underline{x}^{2}=-\|\underline{x}\|^{2}$$ for each
$\underline{x}\in \mathbb{R}^n$.  We have $\mathbb{R}^n\subseteq Cl_{n}$. If $e_1,\ldots, e_n$ is an orthonormal basis for $\mathbb{R}^n$, then $\underline{x}^{2}=-\|\underline{x}\|^{2}$ tells us that $e_i e_j + e_j e_i= -2\delta_{ij}.$ Let $A=\{j_1, \cdots, j_r\}\subset \{1, 2, \cdots, n\}$ and $1\leq j_1< j_2 < \cdots < j_r \leq n$. An arbitrary element of the basis of the Clifford algebra can be written as $e_A=e_{j_1}\cdots e_{j_r}.$ Hence for any element $a\in Cl_{n}$, we have $a=\sum_Aa_Ae_A,$ where $a_A\in \mathbb{R}.$ For $a\in Cl_n$, we will need the anti-involution called
Clifford conjugation:
 $$\bar{a}=\sum_A(-1)^{|A|(|A|+1)/2}a_Ae_A,$$
 where $|A|$ is the cardinality of $A.$
 In particular, we have $\overline{e_{j_1}\cdots e_{j_r}}=(-1)^r e_{j_r}\cdots e_{j_1}$ and $\overline{ab}= \bar{b}\bar{a}$ for $a, b \in Cl_n.$

For each $a=a_0+\cdots +a_{1\cdots n}e_1\cdots e_n\in Cl_n$ the scalar part of $\bar{a}a$ gives the square of the norm of $a,$ namely $a_0^2+\cdots +a_{1\cdots n}^2$\,.

We recall that if $y\in \mathbb{S}^{n-1}\subseteq \mathbb{R}^n$ and $x \in \mathbb{R}^n$, then
$yxy$ gives a reflection of $x$ in the $y$ direction, because
$yxy=yx^{\parallel _y}y+yx^{\perp_y}y=-x^{\parallel _y}+x^{\perp_y}$ where $x^{\parallel _y}$ is the projection of $x$ onto $y$
and $x^{\perp_y}$ is perpendicular to $y$.

\par The Dirac Operator in $\mathbb{R}^n$ is defined to be $$D :=\sum_{j=1}^{n} e_j \frac{\partial}{\partial x_j}.$$
\par Let $\mathcal{M}_k$ denote the space of $Cl_n-$ valued polynomials, homogeneous of degree $k$ and such that if $p_k\in$ $\mathcal{M}_k$ then $Dp_k=0.$ Such a polynomial is called a left monogenic polynomial homogeneous of degree $k$. Note if $h_k\in$ $\mathcal{H}_k,$ the space of $Cl_n-$ valued harmonic polynomials homogeneous of degree $k$, then $Dh_k\in$ $\mathcal{M}$$_{k-1}$. But $Dup_{k-1}(u)=(-n-2k+2)p_{k-1}(u),$
so $$\mathcal{H}_k=\mathcal{M}_k\bigoplus u\mathcal{M}_{k-1}, h_k=p_k+up_{k-1}.$$
This is the so-called Almansi-Fischer decomposition of $\mathcal{H}$$_k$. See \cite{BDS}.

\par Suppose $U$ is a domain in $\mathbb{R}^n$. Consider a function of two variables
$$f: U\times \mathbb{R}^n\longrightarrow Cl_n$$
such that for each $x\in U, f(x,u)$ is a left monogenic polynomial homogeneous of degree $k$ in $u$. Consider the action of the Dirac operator:
$$
D_xf(x,u)\,.
$$
As $Cl_n$ is not commutative then $D_xf(x,u)$ is no longer monogenic in $u$ but it is still harmonic and  homogeneous of degree $k$ in $u$.
So by the Almansi-Fischer decomposition, $D_xf(x,u)=f_{1,k}(x,u)+uf_{2,k-1}(x,u)$ where $f_{1,k}(x,u)$ is a left monogenic polynomial homogeneous of degree $k $ in $u$ and $f_{2,k-1}(x,u)$ is a left monogenic polynomial homogeneous of degree $k-1$ in $u$. Let $P_k$ be the left projection map
 $$P_k:  \mathcal{H}_k\rightarrow \mathcal{M}_k,$$ then
$R_kf(x,u)$ is defined to be $P_kD_xf(x,u)$. The left Rarita-Schwinger equation is defined to be
(see \cite{BSSV})
$$
R_k f(x,u)=0.
$$

\par For an integer $l$, $1\leq l\leq n$ we define the $l$-cylinder $C_l$ to be
the $n$-dimensional manifold $\mathbb{R}^n/\mathbb{Z}^l$, where $\mathbb{Z}^l$ denote
the $l$-dimensional lattice  defined by $\mathbb{Z}^l:= \mathbb{Z}e_1+\cdots+\mathbb{Z}e_l $.
We denote its members $ m_1e_1+\cdots+m_le_l $ for each $m_1,\cdots, m_l\in \mathbb{Z}$
by a bold letter ${\bf m}$. When $l=n$, $C_l$ is the $n$-torus, $T_n.$
For each $l$ the space $\mathbb{R}^n$ is the universal covering of the cylinder $C_l.$
Hence there is a projection map $\pi_l: \mathbb{R}^n \to C_l.$

\par An open subset $U$ of the space $\mathbb{R}^n$ is called $l-$fold periodic if for each $x\in U$ the point $x+{\bf m}\in U$.
So $\pi_l(U)=U'$ is an open subset of the $l$-cylinder $C_l.$

\par Suppose that $U \subset \mathbb{R}^n$ is a $l-$fold periodic open set.
Let $f(x,u)$ be a function defined on $U\times\mathbb{R}^n$ with values in $Cl_n$,
and such that $f$ is a monogenic polynomial homogeneous of degree $k$ in $u$. Then we say that $f(x,u)$ is a $l$-fold periodic function
if for each $x\in U$ we have that $f(x,u)=f(x+{\bf m},u)$.

\par Now if $f: U\times\mathbb{R}^n\to Cl_n$ is a $l-$fold periodic function then the
projection $\pi_l$ induces a well defined function $f': U'\times\mathbb{R}^n\to Cl_n,$
where $f'(x',u)=f(\pi_l^{-1}(x'),u)$ for each $x'=\pi_l(x)\in U'.$ Moreover, any
function $f': U'\times\mathbb{R}^n\to Cl_n$ lifts to a $l-$fold periodic function
$f: U\times\mathbb{R}^n\to Cl_n$, where $U=\pi_l^{-1}(U').$

\par The projection map $\pi_l$ induces a projection of the Rarita-Schwinger
operator $R_k$ to an operator $R_k^{C_l}$ acting on domains on $C_l\times\mathbb{R}^n$ which is
defined by $P_kD'$, where $D'$ is the projection of the Dirac operator $D$.
That is
$$R_k^{C_l}f'(x',u)=P_kD_{x'}'f'(x',u).$$
We call the operator $R_k^{C_l}$ a $l$-cylindrical Rarita-Schwinger type operator and
the solutions of the equation
\begin{equation}\label{rse}
R_k^{C_l}f'(x',u)= 0
\end{equation}
$l$-cylindrical Rarita-Schwinger functions.

As $$I-P_k:  \mathcal{H}_k\rightarrow u\mathcal{M}_{k-1},$$ where $I$ is the identity map, then we can define the remaining Rarita-Schwinger operators
$$Q_k:=(I-P_k)D_x:u\mathcal{M}_{k-1}\to u\mathcal{M}_{k-1}\quad ug(x,u):\to (I-P_k)D_xug(x,u).$$ See\cite{BSSV,LR}.

The remaining Rarita-Schwinger equation is defined to be $(I-P_k)D_xug(x,u)=0$ or $Q_kug(x,u)=0,$ for each $x$ and $(x,u)\in U\times \mathbb{R}^n$, where $U$ is a domain in $\mathbb{R}^n$ and $g(x,u)\in \mathcal{M}_{k-1}.$

Suppose that $g: U\times\mathbb{R}^n\to Cl_n$ is a $l-$fold periodic function then the
projection $\pi_l$ induces a well defined function $g': U'\times\mathbb{R}^n\to Cl_n,$
where $g'(x',u)=g(\pi_l^{-1}(x'),u)$ for each $x'=\pi_l(x)\in U'.$
The projection map $\pi_l$ also induces a projection of the remaining Rarita-Schwinger
operator $Q_k$ to an operator $Q_k^{C_l}$ acting on domains on $C_l\times\mathbb{R}^n$ which is
defined by $(I-P_k)D'$.
That is
$$Q_k^{C_l}ug'(x',u)=P_kD_{x'}'ug'(x',u).$$
We call the operator $Q_k^{C_l}$ a $l$-cylindrical remaining Rarita-Schwinger type operator and
the solutions of the equation
\begin{equation}\label{rse}
Q_k^{C_l}ug'(x',u)= 0
\end{equation}
$l$-cylindrical remaining Rarita-Schwinger functions.

\section{ Fundamental solutions of $R_k^{C_l}$}

Let $U$ a domain in $\mathbb{R}^n$. We recall the fundamental solution of the Rarita-Schwinger operator $R_k$ in
$\mathbb{R}^n$, see \cite{DLRV}:
$$
E_k(x,u,v)=\di\frac{1}{\omega_n c_k}\di\frac{x}{\|x\|^n}Z_k(\frac{xux}{\|x\|^2},v)\,,
$$
where $c_k=\di\frac{n-2}{n+2k-2}, \omega_n$ is the surface area of the unit sphere in $\mathbb{R}^n,$
$$
Z_k(u, v):= \sum_{\sigma}{P_{\sigma}(u)V_{\sigma}(v)v}\,,
$$
where $P_{\sigma}(u)= \di\frac{1}{k!}\sum_{\sigma}{(u_{i_1}-u_1e_1^{-1}e_{i_1})\ldots (u_{i_k}-u_1e_1^{-1}e_{i_k})}$\,,
$V_{\sigma}(v) = \frac{\partial^k G(v)}{\partial v_{i_2}^{k_2} \ldots \partial v_{i_k}^{k_n}}, $
 $k_2+\ldots +k_n=k, i_k\in\{2,\cdots,n\}$ and the summation is taken over all
permutations of the monomials without repetition. See \cite{BDS}.

Now we construct the following functions
\begin{equation}\label{cot1}
\cot_{l,k}(x,u,v)=\di\sum_{(m_1,\cdots,m_l)\in\mathbb{Z}^l}E_k(x+m_1e_1+\cdots+m_le_l,u,v)~, ~~ \mbox{for}~~ 1\leq l\leq n-2.
\end{equation}

These functions are defined on the $l$-fold periodic domain $\mathbb{R}^n/\mathbb{Z}^l$ for fixed $u$ and
$v$ in $\mathbb{R}^n$ and  are $Cl_n$-valued. It is easy to see that they are
$l$-fold periodic functions.

Now, we will prove the locally uniform convergence of the series $\cot_{l,k}(x,u,v)$. First, we give a detailed proof
for the locally normal convergence of the series $\di\sum_{\bf m \in\mathbb{Z}^l}{G(x + \bf m)}$,
where $G(x + {\bf m}) =\di\frac{x + {\bf m}}{||x + {\bf m}||^n}$.
We will need the following proposition whose proof can be found on page 42 in \cite{K}.
\begin{proposition}\label{1}
Let $N_0$ be the set of non-negative integers and ${\cal A}_{k+1}$ be the space of paravectors $z = x_0 + {\bf x}$
with $Sc(z)= x_0$ and $Vec(z)= {\bf x} \in \mathbb{R}^k$. Let $s \in \{1, \dots, k\}$. For all multi-indices
${\bf \alpha} \in \mathbb{N}^{k+1}_0$, $|{\bf\alpha}| \geq 1$, ${\bf \alpha} = (0, \alpha_1, \dots, \alpha_k)$, the following estimate holds
for all $z \in {\cal A}_{k+1}\setminus \{0\} $:
$$
||\frac{\partial^{|{\bf \alpha}|}}{\partial z^{\bf \alpha}} q_{\bf 0}^{(s)} (z)|| \leq
\frac{(k+1-s)(k+2-s)\ldots(k +|{\bf \alpha}|-s)}{||z||^{k+|{\bf \alpha}|+1-s}}\,,
$$
where ~~
$q_{\bf 0}(z):= \displaystyle \frac{\bar{z}}{||z||^{k+1}} $\,
$\displaystyle \frac{\partial^{|{\bf \alpha}|}}{\partial z^{\bf \alpha}} q_{\bf 0}(z):=
\frac{\partial^{\alpha_1 + \ldots + \alpha_k}}{\partial x_1^{\alpha_1}\ldots \partial x_k^{\alpha_k}} q_{\bf 0}(z)$\,
~and~ $q_{\bf 0}^{(s)}(z)$ is the kernel of the $s$th-power of the Dirac operator.
\end{proposition}
\begin{remark}
In the vector formalism in  $\mathbb{R}^n $, $ q_{\bf 0}({\bf x}): = \displaystyle-\frac{{\bf x}}{||{\bf x}||^n} $.
\end{remark}

From now on we are working only in the case $s=1$.

Now we have a particular case of Proposition 2.2 appearing in \cite{K}.

\begin{proposition}\label{2}
Let $p\in \mathbb{N}$ with $1 \leq p \leq n-2$. Let $\mathbb{Z}^p$ be the $p$-dimensional lattice. Then the series
\begin{equation}\label{series1}
\sum_{{\bf m}\in \mathbb{Z}^p}{q_{\bf 0}(x + {\bf m})}
\end{equation}
converges normally in $\mathbb{R}^n \setminus \mathbb{Z}^p $.
\end{proposition}
\proof
We consider an arbitrary compact subset $ K\subset \mathbb{R}^n $ and a real number $R > 0$ such that the ball
$\bar{B}(0, R)$ covers $K$ completely. Let $x \in \bar{B}(0, R), x= (x_1, \dots, x_n) $. WLOG we will study the convergence of the series
summing only over those lattice ${\bf m}$ that satisfy $|| {\bf m}|| > n R \geq ||x||$.

The function $q_{\bf 0}(x + {\bf m})$ is left monogenic in $0 \leq ||x|| < n R$. Hence it is real analytic in
$\bar{B}(0, n R) $ and therefore can be represented in the interior of this ball by its Taylor series, i.e.,
$$
q_{\bf 0}(x + {\bf m}) =
\sum_{\nu=0}^{\infty}\Big( {\sum_{l_1 + \ldots + l_n = \nu} {\frac{1}{\bf l!}x_1^{l_1} \ldots x_n^{l_n} q_{{\bf l}}({\bf m})}}\Big),
$$
where ${\bf l}=(l_1,\cdots,l_n)$ and $q_{{\bf l}}({\bf m})=
\di\frac{\partial ^{|\bf l|}}{\partial{\bf m}^{\bf l}}q_0({\bf m}).$

Using Proposition \ref{1} and observing that $||x||^{\nu} \leq R^{\nu}$, we obtain:
\begin{eqnarray*}
||q_{\bf 0}(x + {\bf m})|| & \leq &
\sum_{\nu=0}^{\infty} \Big({\sum_{l_1 + \ldots + l_n = \nu} {\frac{1}{\bf l!} ||x||^{\nu} ||q_{{\bf l}}({\bf m})||}}\Big) \\
& \leq & \sum_{\nu=0}^{\infty} \Big({\sum_{l_1 + \ldots + l_n = \nu} {\frac{1}{l_1! \ldots l_n!}\, R^{\nu}\, \nu!\,
\prod_{\gamma = 1}^{n-2}(\nu + \gamma) \frac{1}{||{\bf m}||^{n-1 + \nu}}}}\Big)\,.
\end{eqnarray*}
Using the multinomial formula at
$$
\sum_{l_1 + \ldots + l_n = \nu} \frac{\nu!}{l_1! \ldots l_n!} 1^{l_1}\cdots 1^{l_n} = n^{\nu} \,,
$$
we can write the former estimate as
$$
||q_{\bf 0}(x + {\bf m})|| \leq
\sum_{\nu=0}^{\infty}{ \prod_{\gamma = 1}^{n-2}(\nu + \gamma)\Big(\frac{n R}{||{\bf m}||}\Big)^{\nu}\frac{1}{||{\bf m}||^{n-1}}}\,.
$$
Since the series $\displaystyle\sum_{k=-s}^{\infty}{r^{k+s}}$ converges absolutely for $|r| < 1$ to $\displaystyle \frac{1}{1-r}$,
we can consider its $s$th-derivative which also converges absolutely for $|r| < 1$ to
$\displaystyle \frac{d^s}{dr^s}(\frac{1}{1-r})$. So we obtain
$\displaystyle \sum_{k=0}^{\infty}{(k+1)(k+2)\cdots(k+s)r^k} = \frac{s!}{(1-r)^{s+1}}\,.$
Taking this into account and observing that $\displaystyle \frac{nR}{||{\bf m}||} < 1$, we get
$$
\sum_{\nu=0}^{\infty}{\prod_{\gamma = 1}^{n-2}(\nu + \gamma)\Big(\frac{n R}{||{\bf m}||}\Big)^{\nu}} =
\frac{(n-2)!}{(1 - \frac{nR}{||{\bf m}||})^{n-1}}\,\,.
$$
Therefore we have
$$
||q_{\bf 0}(x + {\bf m})|| \leq \frac{(n-2)!}{{(1 - \frac{nR}{||{\bf m}||})^{n-1}}} \cdot \frac{1}{||{\bf m}||^{n-1}}\,.
$$
Due to Eisenstein's Lemma (see \cite{E}) a series of the form
$$
\sum_{(m_1, \ldots, m_p)\in \mathbb{Z}^p\setminus \{0\}}{||m_1\omega_1 + \ldots + m_p\omega_p||^{-(p+\alpha)}}
$$
is convergent if and only if $\alpha \geq 1$, where $\omega_i, i= 1, \ldots, p,$ are $\mathbb{R}$-linear independent paravectors in ${\cal A}_{n}.$

In our case, $1\leq p \leq n-2$ which implies \- $\displaystyle \frac{A}{||{\bf m}||^{n-1}}= \frac{A}{||{\bf m}||^{p+\alpha}}$
\-\- for $\alpha \geq 1$ and \-\- $A = \di\frac{(n-2)!}{{\Big(1 - \frac{nR}{||{\bf m}||}\Big)^{n-1}}}$,
and therefore observing that also $e_j, j= 1, \ldots, p,$ are
$\mathbb{R}$-linear independent paravectors in ${\cal A}_{n}$ and using the comparison test, we obtain that the series
(\ref{series1}) converges normally in $\mathbb{R}^n \setminus \mathbb{Z}^p $.  $\blacksquare$

Returning to the series defined by (\ref{cot1})
\begin{eqnarray*}
\cot_{l,k}(x,u,v)& = &\di\sum_{{\bf m}\in \mathbb{Z}^l}E_k(x+m_1e_1+\cdots+m_le_l,u,v) \\
& = & \sum_{{\bf m}\in \mathbb{Z}^l}{G( x + {\bf m})} Z_k\Big(\frac{(x + {\bf m})u (x + {\bf m})}{\|x + {\bf m}\|^2}, v \Big), 1\leq l\leq n-2,
\end{eqnarray*}
we observe that
$Z_k(\di\frac{(x+{\bf m})u(x+{\bf m})}{\|x+{\bf m}\|^2},v)$ is a bounded function on a bounded domain in $\mathbb{R}^n$, because its first variable, $\di\frac{(x + {\bf m})u (x + {\bf m})}{\|x + {\bf m}\|^2}$, is a reflection
in the direction $\displaystyle\frac{x + {\bf m}}{||x + {\bf m}||}$ for each ${\bf m}$,
and hence is a linear transformation which is a continuous function. On the other hand, we would get
with respect to the second variable, bounded homogeneous functions of degree $k$.

Consequently, applying Proposition \ref{2} the series (\ref{cot1}) is a uniformly convergent
series and represents a kernel for the Rarita-Schwinger operators under translations by
${\bf m}\in \mathbb{Z}^{l}$, with $1\leq l\leq n-2 $.

Now, we want to define the $(n-1)$-fold periodic cotangent. In order to do that, we decompose as in \cite{K} the lattice $\mathbb{Z}^{l}$
into three parts: the origin $\{{\bf 0}\}$ and a positive and a negative parts. The last two parts are equal and disjoint:
\begin{eqnarray*}
\Lambda_l & = &\{m_1e_1: m_1\in \mathbb{N}\}\cup\{m_1e_1+m_2e_2: m_1, m_2\in\mathbb{Z}, m_2>0\} \\
&& \cup\cdots\cup\{m_1e_1+\cdots+m_le_l: m_1,\cdots, m_l\in\mathbb{Z}, m_l>0\} \\
-\Lambda_l & = & (\mathbb{Z}^{l}\setminus \{0\})\setminus \Lambda_l.
\end{eqnarray*}
For $l=n-1,$ we define
\begin{equation}\label{cot11}
\cot_{n-1,k}(x,u,v)=E_k(x,u,v)+\di\sum_{{\bf m}\in \Lambda_{n-1}}[E_k(x+{\bf m},u,v)+ E_k(x-{\bf m},u,v)].
\end{equation}

To show the uniform convergence of the above series, we need the following proposition
which is a special case in \cite{K}.

\begin{proposition}\label{4}
Let $\mathbb{Z}^{n-1}$ be the $(n-1)$-dimensional lattice. Then the series
\begin{equation}\label{series11}
q_{\bf 0}(x) + \sum_{{\bf m}\in \mathbb{Z}^{n-1}\setminus \{0\}}{\big(q_{\bf 0}(x + {\bf m})- q_{\bf 0}({\bf m})\big)}\,,
\end{equation}
converges normally in $\mathbb{R}^n \setminus \mathbb{Z}^{n-1}$.
\end{proposition}
\proof Following the proof of Proposition \ref{2}, again we suppose that $x \in \bar {B}(0, R)$ and
consider only those lattice points with
$|| {\bf m}|| > n R \geq ||x||$. We consider the function
$q_{\bf 0}(x + {\bf m}) - q_{\bf 0}({\bf m})$ and expand it into a Taylor series in $B(0, R')$
where $0 < R' < R$. So we have
$$
q_{\bf 0}(x + {\bf m}) - q_{\bf 0}({\bf m}) =
\sum_{\nu=1}^{\infty} \Big({\sum_{l_1 + \ldots + l_{n} = \nu}{\frac{1}{\bf l!}x_1^{l_1} \ldots x_n^{l_{n}} q_{{\bf l}}({\bf m})}}\Big).
$$
Using similar arguments to those of proof of Proposition \ref{2}, we obtain
\begin{eqnarray*}
||q_{\bf 0}(x + {\bf m}) - q_{\bf 0}({\bf m})|| &\leq &
\sum_{\nu=1}^{\infty}{\prod_{\gamma = 1}^{n-2}(\nu + \gamma)\Big(\frac{n R}{||{\bf m}||}\Big)^{\nu}\frac{1}{||{\bf m}||^{n-1}}} \\
& = &
\sum_{\nu=0}^{\infty}{\prod_{\gamma = 2}^{n-1}(\nu + \gamma)\Big(\frac{n R}{||{\bf m}||}\Big)^{\nu}\frac{nR}{||{\bf m}||^n}}\,.
\end{eqnarray*}
Since
$$
||q_{\bf 0}(x + {\bf m}) - q_{\bf 0}({\bf m})|| \leq \frac{KR}{||{\bf m}||^n} = \frac{KR}{||{\bf m}||^{p+1}}\,,~~ p=n-1,
$$
where $K$ is a positive real constant and applying the Eisenstein's Lemma, the series (\ref{series11})
turns into a normally convergent series
in $\mathbb{R}^n \setminus \mathbb{Z}^{n-1}$. $\blacksquare$
\begin{remark}
Because $\Lambda_{n-1}$ and $-\Lambda_{n-1}$ belong to
$\mathbb{Z}^{n-1}\setminus \{0\}$, we have also the normal convergence of the
series
$$
\sum_{{\bf m}\in \Lambda_{n-1}}{\big(q_{\bf 0}(x + {\bf m})- q_{\bf 0}({\bf m})\big)}, ~~
\sum_{{\bf m}\in -\Lambda_{n-1}}{\big(q_{\bf 0}(x + {\bf m})- q_{\bf 0}({\bf m})\big)}
$$
\end{remark}
Now we consider the expression
$$
G( x + {\bf m}) Z_k\Big(\frac{(x + {\bf m})u (x + {\bf m})}{\|x + {\bf m}\|^2}, v \Big)
+G( x - {\bf m}) Z_k\Big(\frac{(x - {\bf m})u (x - {\bf m})}{\|x - {\bf m}\|^2}, v \Big)\,,
$$
where the function $Z_k(\cdot, \cdot)$ has the hypotheses stated formerly. Rewriting
the former expression as
\begin{multline*}
(G( x + {\bf m})+ G( x - {\bf m}))Z_k\Big(\frac{(x + {\bf m})u (x + {\bf m})}{\|x + {\bf m}\|^2}, v \Big)\\
+ (-G( x - {\bf m}) - G({\bf m}))Z_k\Big(\frac{(x + {\bf m})u (x + {\bf m})}{\|x + {\bf m}\|^2}, v \Big)\\
+ (G( x - {\bf m}) + G({\bf m})) Z_k\Big(\frac{(x - {\bf m})u (x - {\bf m})}{\|x - {\bf m}\|^2}, v \Big) \\
+ G({\bf m})Z_k\Big(\frac{(x + {\bf m})u (x + {\bf m})}{\|x + {\bf m}\|^2}, v \Big)
- G({\bf m})Z_k\Big(\frac{(x - {\bf m})u (x - {\bf m})}{\|x - {\bf m}\|^2}, v \Big)
\end{multline*}
and taking the sum over $\mathbb{Z}^{n-1}\setminus \{0\}$ we obtain
\begin{multline}\label{caseb}
\sum_{{\bf m}\in \mathbb{Z}^{n-1}\setminus \{0\}}\Big(G( x + {\bf m}) Z_k\Big(\frac{(x + {\bf m})u (x + {\bf m})}{\|x + {\bf m}\|^2}, v \Big)
+G( x - {\bf m}) Z_k\Big(\frac{(x - {\bf m})u (x - {\bf m})}{\|x - {\bf m}\|^2}, v \Big)\Big) = \\
\sum_{{\bf m}\in \mathbb{Z}^{n-1}\setminus \{0\}}(G( x + {\bf m})+ G( x - {\bf m}))
Z_k\Big(\frac{(x + {\bf m})u (x + {\bf m})}{\|x + {\bf m}\|^2}, v \Big)\\
+ \sum_{{\bf m}\in \mathbb{Z}^{n-1}\setminus \{0\}}(-G( x - {\bf m}) - G({\bf m}))
Z_k\Big(\frac{(x + {\bf m})u (x + {\bf m})}{\|x + {\bf m}\|^2}, v \Big)\\
+ \sum_{{\bf m}\in \mathbb{Z}^{n-1}\setminus \{0\}}(G( x - {\bf m}) + G({\bf m}))
Z_k\Big(\frac{(x - {\bf m})u (x - {\bf m})}{\|x - {\bf m}\|^2}, v \Big) \\
+ \sum_{{\bf m}\in \mathbb{Z}^{n-1}\setminus \{0\}}\Big(G({\bf m})Z_k\Big(\frac{(x + {\bf m})u (x + {\bf m})}{\|x + {\bf m}\|^2}, v \Big)
- G({\bf m})Z_k\Big(\frac{(x - {\bf m})u (x - {\bf m})}{\|x - {\bf m}\|^2}, v \Big)\Big).
\end{multline}

The last sum in (\ref{caseb}) vanishes because the terms
$G({\bf m})Z_k\Big(\frac{(x + {\bf m})u (x + {\bf m})}{\|x + {\bf m}\|^2}, v \Big)$
for each ${\bf m} \in \Lambda_{n-1}$ and ${\bf m} \in -\Lambda_{n-1}$ cancel with the terms
$G({\bf m})Z_k\Big(\frac{(x - {\bf m})u (x - {\bf m})}{\|x - {\bf m}\|^2}, v \Big)$
for each ${\bf m} \in -\Lambda_{n-1}$ and ${\bf m} \in \Lambda_{n-1}$, respectively.
On the other side, since
\begin{eqnarray*}
||G( x + {\bf m}) + G( x - {\bf m})||
\leq  || G( x + {\bf m}) - G({\bf m})|| + || G( x + \overline{{\bf m}}) - G(\overline{{\bf m}})||.
\end{eqnarray*}
and observing that $Z_k(\cdot, \cdot)$ is bounded, we can apply
Proposition \ref{4} to obtain the normal convergence for the series (\ref{caseb}).
Consequently the series defined by (\ref{cot11}) is uniformly convergent and it is a
kernel for the Rarita-Schwinger type operator
under translations by  ${\bf m}\in \Lambda_{n-1}$.

For $x, y\in \mathbb{R}^n\setminus \mathbb{Z}^l$, $1\leq l \leq n-1$,
the functions $\cot_{l,k}(x-y,u,v)$ induce functions $$\cot_{l,k}'(x',y',u,v)=\cot_{l,k}(\pi_l^{-1}(x')-\pi_l^{-1}(y'),u,v).$$
These functions are defined on $(C_l\times C_l)\setminus diag(C_l\times C_l)$ for each fixed $u,v\in \mathbb{R}^n,$ where $diag(C_l\times C_l)=\{(x',x'):x'\in C_l\}$ and they are $l$-cylindrical Rarita-Schwinger functions, i.e,
the equation $R_k^{C_l}\cot_{l,k}'(x',y',u,v)=0$ is satisfied. Furthermore for each $l$ they represent a kernel for the
operator $R_k^{C_l}$.

\section{Some integral formulas on cylinders}

\begin{definition} For any $Cl_n$-valued polynomials $P, Q$, the inner product
$(P(u), Q(u))_u$ with respect to $u$ is given by $$(P(u), Q(u))_u=\di\int_{\mathbb{S}^{n-1}}P(u)Q(u)dS(u),$$
where $\mathbb{S}^{n-1}$ is the unit sphere in $\mathbb{R}^n$.\end{definition} See\cite{DLRV}.

\par For any $p_k \in \mathcal{M}_k,$ one obtains $$p_k(u)=(Z_k(u,v), p_k(v))_v=\int_{\mathbb{S}^{n-1}}Z_k(u,v)p_k(v)dS(v).$$
See \cite{BDS}.
\begin{theorem}\cite{DLRV} (Stokes' Theorem for $R_k$)  Let $\Omega$ and $\Omega'$ be domains in $\mathbb{R}^n$
and suppose the closure of $\Omega$ lies in $\Omega'$. Further suppose the closure of $\Omega$ is
compact and $\partial\Omega$ is piecewise smooth. Then for $f, g \in C^1(\Omega',\mathcal{M}_k)$, we have
$$\begin{array}{ll}
\di\int_\Omega[(g(x,u)R_k, f(x,u))_u+(g(x,u), R_kf(x,u))_u]dx^n\\
\\
=\di\int_{\partial\Omega}\left(g(x,u), P_kd\sigma_xf(x,u)\right)_u,
\end{array}$$ where $dx^n = dx_1\wedge\cdots \wedge dx_n$, $d\sigma_x=\di\sum_{j=1}^n(-1)^{j-1}e_jd\hat{x_j},$
 and $d\hat{x_j}=dx_1\wedge\cdots dx_{j-1}\wedge dx_{j+1}\cdots \wedge dx_n$.\end{theorem}
\begin{theorem} \cite{DLRV} (Borel-Pompeiu Theorem) Let $\Omega'$ and $\Omega$ be as in Theorem 1.
Then for $f \in C^1(\Omega',\mathcal{M}_k)$
$$\begin{array}{ll}
f(y,v)=\di\int_{\partial\Omega}\left(E_k(x-y,u,v), P_kd\sigma_xf(x,u)\right)_u-\di\int_\Omega(E_k(x-y,u,v),R_kf(x,u))_udx^n.
\end{array}$$\end{theorem}
\par Now by Stokes' Theorem and Borel-Pompeiu Theorem for the Rarita-Schwinger operator $R_k$ in $\mathbb{R}^n$, we may easily obtain:

\begin{theorem} \quad Let $V$ be a bounded domain in $\mathbb{R}^n$ and $\overline{V}$ be
the closure of $V$. For each $x\in \overline{V},$ the shifted lattice $x+\mathbb{Z}^l$
intersected with $V$ only contains the points $x.$ Suppose that the boundary of $V$, $\partial {V}$, is piecewise
smooth and $\overline{V}$ is compact. Further suppose $f(x,u): \overline{V}\times \mathbb{R}^n\to Cl_n$ is a
monogenic homogeneous polynomial of degree $k$ in $u$ and with respect to $x$ is $C^1$. Then for $1\leq l\leq n-1$ and each $y\in V,$
$$
f(y,v)=\di\int_{\partial V}\left(\cot_{l,k}(x-y,u,v), P_kd\sigma_xf(x,u)\right)_u-\di\int_V(\cot_{l,k}(x-y,u,v),R_kf(x,u))_udx^n.
$$
\end{theorem}

\proof.  When ${\bf m}=0,$ this is the Borel-Pompeiu Formula given by Theorem 2
\begin{equation}\label{eqn1}
f(y,v)=\di\int_{\partial V}(E_k(x-y,u,v), P_kd\sigma_xf(x,u))_u-\di\int_V(E_k(x-y,u,v),R_kf(x,u))_udx^n.
\end{equation}
\par When ${\bf m}\neq 0,$ for each ${\bf m}\in \mathbb{Z}^l$, by the hypothesis
that the shifted lattice $x+\mathbb{Z}^l$ intersected with $V$ only contains the
points $x$, we have $x+{\bf m} \notin \overline{V}$ but $y\in \overline{V}$, so there
is no singularity in $E_k(x-y+{\bf m},u,v)=E_k(x+{\bf m}-y,u,v).$ Hence by Stokes' Theorem we obtain
\begin{eqnarray}\label{eqn2}\di\int_{\partial V}\left(E_k(x-y+{\bf m},u,v), P_kd\sigma_xf(x,u)\right)_u-\di\int_V(E_k(x-y+{\bf m},u,v),R_kf(x,u))_udx^n\nonumber\\
\nonumber\\
=\di\int_V(E_k(x-y+{\bf m},u,v)R_k,f(x,u))_udx^n=0.
\end{eqnarray}
Now for all ${\bf m}\in \mathbb{Z}^l$, by adding the equations (\ref{eqn1}) and (\ref{eqn2}), we obtain
$$
\begin{array}{ll}
f(y,v)=\di\int_{\partial V}\left(\di\sum_{{\bf m}\in\mathbb{Z}^l}E_k(x-y+{\bf m},u,v), P_kd\sigma_xf(x,u)\right)_u\\
\\
\qquad\qquad-\di\int_V(\di\sum_{{\bf m}\in\mathbb{Z}^l}E_k(x-y+{\bf m},u,v),R_kf(x,u))_udx^n.\qquad \blacksquare
\end{array}
$$

\par The case we are most interested in here is the one that $V$ is $l-$fold periodic and $f$ is $l-$fold periodic. If we use
the projection map $\pi_l$ at this moment, we obtain the Borel-Pompeiu Theorem on the cylinder $C_l.$

\begin{theorem}\label{BPT}\quad (Borel-Pompeiu Theorem for $R_k^{C_l}$)
Suppose $V'$ is a domain in $C_l$ with compact closure and smooth boundary. Suppose $f(x,v)$ is defined as in Theorem 3.
Then for $1\leq l\leq n-1$ and each $y'\in V',$
$$\begin{array}{ll}
f'(y',v)=\di\int_{\partial V'}\left(\cot_{l,k}'(x',y',u,v), P_kd\sigma'_{x'}f'(x',u)\right)_u\\
\\
\qquad \quad -\di\int_{V'}(\cot_{l,k}'(x',y',u,v),R_k^{C_l}f'(x',v))_vd\mu(x'),
\end{array}$$
where $x'=\pi_l(x), d\sigma'_{x'}=\partial_x\pi_ld\sigma_{x},$ $\partial_x\pi_l$ is the derivative of $\pi_l$ at $x$,
 $\mu$ is the projection of Lebesgue measure on $\mathbb{R}^n$ onto $C_l,$ $d\sigma_{x}= n(x)d\sigma(x)$, $n(x)$ is the unit exterior normal vector at $x = \pi_l^{-1}(x')$ and $d\sigma(x)$ is the surface measure
element.
\end{theorem}
\begin {theorem}\quad (Cauchy integral formula)
Suppose the hypotheses as in Theorem \ref{BPT} and $f'(x',u)$ is annihilated by
the operator $R_k^{C_l},$ then for $1\leq l\leq n-1$ and each $y'\in V',$
we have $$f'(y',u)=\di\int_{\partial V'}\left(\cot_{l,k}'(x',y',u,v), P_kd\sigma'_{x'}f'(x',u)\right)_u.$$
\end{theorem}

\par If the function given in the Borel-Pompeiu Theorem has its support with respect to $x'$
in $V'\subset C_l,$ then by the same theorem, we have the following:
\begin{theorem}
\quad $\di\int_{V'}-(\cot_{l,k}'(x',y',u,v),R_k^{C_l}\psi (x',u))_udx'^n=\psi(y',v),$ for $\psi\in C^{\infty}(C_l).$
\end{theorem}

\par We also can introduce a Cauchy transform for the Rarita-Schwinger operator $R_k^{C_l}:$
\begin{definition} For a subdomain $V'$ of cylinder $C_l$ and a function $f': V' \longrightarrow Cl_n,$ the Cauchy
transform of $f'$ is formally defined to be
$$
(Tf')(y',v)=-\int_{V'} \left(\cot_{l,k}'(x',y',u,v), f'(x',u)\right)_udx'^n,\qquad y'\in V'.
$$
\end{definition}
\par Consequently, one may obtain a right inverse for $R_k^{C_l}$
\begin{theorem}
\quad $R_k^{C_l}\di\int_{V'}(\cot_{l,k}'(x',y',u,v),\psi (x',u))_udx'^n=\psi(y',v),$ for $\psi\in C^{\infty}(C_l).$
\end{theorem}

\section{Conformally inequivalent spinor bundles on $C_l$}

We shall now show a construction of a number of conformally inequivalent spinor bundles on $C_l$.
In the previous sections the spinor bundle over $C_l$ is chosen to be the trivial one $C_l \times Cl_n$,
however we can construct $2^{l}$ spinor bundles on $C_l$.

Different spin structures on a
spin manifold $M$ are detected by the number of different homomorphisms from the
fundamental group $\Pi_1(M)$ to the group $\mathbb{Z}_2=\{0,1\}$. In our case we have
$\Pi_1(C_l)= \mathbb{Z}^l$. Because there are two homomorphisms of $\mathbb{Z}$ to $\mathbb{Z}_2$,
we have $2^{l}$ distinct spin structures on $C_l$, see \cite{MP, KR1}.

The following construction is for some of spinor bundles over $C_l$ but all the others are constructed
similarly. First let $p$ be an integer in the set $\{1, \dots, l\}$ and consider the lattice
$\mathbb{Z}^p:= \mathbb{Z}e_1+\cdots+\mathbb{Z}e_p$. We also consider the lattice
$\mathbb{Z}^{l-p}:= \mathbb{Z}e_{p+1}+\cdots+\mathbb{Z}e_l$. In this case
$\mathbb{Z}^l= \{{\bf m}+ {\bf n}: {\bf m}\in \mathbb{Z}^p ~~ \mbox{and} ~~  {\bf n}\in \mathbb{Z}^{l-p}\}$.
Suppose that ${\bf m} = m_1e_1+\cdots+ m_pe_p$. Let us make the identification $(x, X)$ with
$(x + {\bf m} + {\bf n}, (-1)^{m_1 + \cdots + m_p}X)$ where $ x\in \mathbb{R}^n $ and $X\in Cl_n$.
This identification gives rise to a spinor bundle $E^p$ over $C_l$.

The Rarita-Schwinger operator over $\mathbb{R}^n$ induces a Rarita-Schwinger operator acting on sections
of the bundles $E^p$ over $C_l$. We will denote this operator by $R_{k,p}^{C_l}$.
The projection $\pi_l$ maps $U\subset \mathbb{R}^n $ to a domain $U'\subset C_l$.
Now if $f: U\times\mathbb{R}^n\to Cl_n$ is a $l-$fold periodic function then the
projection $\pi_l$ induces a well defined function $f'_p: U'\times\mathbb{R}^n\to E^p,$
where $f'_p(x'+ {\bf m} + {\bf n},u)= (-1)^{m_1 + \cdots + m_p}f(\pi_l^{-1}(x')+{\bf m} + {\bf n} ,u)$ for each
$x'+{\bf m} + {\bf n}=\pi_l(x)+{\bf m} + {\bf n}\in U'.$
If $R_{k,p}^{C_l}(f'_p)=0$ then $f'_p$ is called an $E^p$ left Rarita-Schwinger section.
Moreover, any $E^p$ left Rarita-Schwinger section
$f'_p: U'\times\mathbb{R}^n\to E^p$ lifts to a $l-$fold periodic function
$f: U\times\mathbb{R}^n\to Cl_n$, where $U=\pi_l^{-1}(U')$, and $R_k^{C_l}f=0$.

By considering the series
$$
\cot_{l,k}(x,u,v)=\di\sum_{{\bf m}\in\mathbb{Z}^l}E_k(x + {\bf m} ,u,v)~, ~~ \mbox{for}~~ 1\leq l\leq n-2
$$
which converge normally on $\mathbb{R}^n \setminus \mathbb{Z}^l $, we can obtain the kernel (see section 3)
$$
\cot_{l,k}(x,y,u,v)=\di\sum_{{\bf m}\in\mathbb{Z}^l}E_k(x-y + {\bf m} ,u,v)~, ~~ \mbox{for}~~ 1\leq l\leq n-2.
$$
Applying the projection map $ \pi_l:\mathbb{R}^n \to C_l $ to these kernels induce kernels
$$
\cot'_{l,k}(x',y',u,v)
$$
defined on $(C_l \times C_l) \setminus \mbox{diag}(C_l\times C_l)$, where
$\mbox{diag}(C_l\times C_l)= \{(x', x'): x'\in C_l\}$. We can adapt these functions as follows. For
$1\leq l \leq n-2$ we define
$$
\cot_{l,k,p}(x,u,v)=\di\sum_{{\bf m}\in\mathbb{Z}^p, {\bf n}\in\mathbb{Z}^{l-p}}(-1)^{m_1 + \cdots + m_p}E_k(x + {\bf m} + {\bf n} ,u,v).
$$
These are well defined functions on $\mathbb{R}^n \setminus \mathbb{Z}^l $. Therefore we obtain from these functions
the cotangent kernels
$$
\cot_{l,k,p}(x,y,u,v)=\di\sum_{{\bf m}\in\mathbb{Z}^p, {\bf n}\in\mathbb{Z}^{l-p}}(-1)^{m_1 + \cdots + m_p}E_k(x-y + {\bf m}+ {\bf n} ,u,v).
$$
Again applying the projection map $ \pi_l$ these kernels give rise to the kernels
$$
\cot'_{l,k,p}(x',y',u,v).
$$
In the case $l=n-1$, by considering the series
$$
\cot_{n-1,k}(x,u,v)=E_k(x,u,v)+\di\sum_{{\bf m}\in \mathbb{Z}^{n-1}}[E_k(x+{\bf m},u,v)+E_k(x-{\bf m},u,v)]
$$
we obtain the kernel
\begin{multline*}
\cot_{n-1,k}(x,y,u,v)=E_k(x-y,u,v)+\di\sum_{{\bf m}\in \mathbb{Z}^{n-1}}[E_k(x-y + {\bf m},u,v)\\
+E_k(x-y - {\bf m},u,v)]
\end{multline*}
which in turn using the projection map induces kernels $\cot'_{n-1,k}(x',y',u,v)$.
Defining
\begin{multline*}
\cot_{n-1,k,p}(x,u,v)= E_k(x + {\bf m} + {\bf n} ,u,v) + \\
\sum_{{\bf m}\in\mathbb{Z}^p, {\bf n}\in\mathbb{Z}^{n-1-p}}(-1)^{m_1 + \cdots + m_p}\Big[E_k(x + {\bf m} + {\bf n} ,u,v)
+ E_k(x - {\bf m} - {\bf n} ,u,v)\Big]
\end{multline*}
we obtain the cotangent kernels
\begin{multline*}
\cot_{n-1,k,p}(x,y,u,v)= E_k(x -y + {\bf m} + {\bf n} ,u,v) + \\
\sum_{{\bf m}\in\mathbb{Z}^p, {\bf n}\in\mathbb{Z}^{n-1-p}}(-1)^{m_1 + \cdots + m_p}\Big[E_k(x-y + {\bf m} + {\bf n} ,u,v)
+ E_k(x -y - {\bf m} - {\bf n} ,u,v)\Big]
\end{multline*}
and by $\pi_l$ the kernels $\cot'_{n-1,k,p}(x',y',u,v)$.

\section{ Remaining Rarita-Schwinger type operators on cylinders }

Consider the fundamental solution of the remaining Rarita-Schwinger operator $Q_k$ in
$\mathbb{R}^n$, see \cite{LR}:
$$H_k(x,u,v):=\di\frac{-1}{\omega_n c_k}u\di\frac{x}{\|x\|^n}Z_{k-1}(\di\frac{xux}{\|x\|^2},v)v,$$where $c_k=\di\frac{n-2}{n-2+2k}.$

Now we construct functions
\begin{multline*}
Cot_{l,k}(x,u,v)=\di\sum_{{\bf m}\in \mathbb{Z}^l}H_k(x+m_1e_1+\cdots+m_le_l,u,v) \\
= \sum_{{\bf m}\in \mathbb{Z}^l}u{G( x + {\bf m})} Z_{k-1}\Big(\frac{(x + {\bf m})u (x + {\bf m})}{\|x + {\bf m}\|^2}, v \Big)v, 1\leq l\leq n-2.
\end{multline*}

These functions are defined on the $l$-fold periodic domain $\mathbb{R}^n/\mathbb{Z}^l$ for fixed $u$ and
$v$ in $\mathbb{R}^n$ and  are $Cl_n$-valued. We can observe that they are
$l$-fold periodic functions. Using the similar arguments in Section $3,$ we may easily obtain that the series $Cot_{l,k}(x,u,v)$ is normally convergent over $\mathbb{R}^n \setminus \mathbb{Z}^l $.

For $l=n-1,$ we define
$$
Cot_{n-1,k}(x,u,v)=H_k(x,u,v)+\di\sum_{{\bf m}\in \mathbb{Z}^{n-1}}[H_k(x+{\bf m},u,v)+H_k(x-{\bf m},u,v)].
$$
We can establish that the previous series converges uniformally following the proof in Section 3 for the case $l=n-1.$

For $x, y\in \mathbb{R}^n\setminus \mathbb{Z}^l$, $1\leq l \leq n-1$,
the functions $Cot_{l,k}(x-y,u,v)$ induce functions $$Cot_{l,k}'(x',y',u,v)=Cot_{l,k}(\pi_l^{-1}(x')-\pi_l^{-1}(y'),u,v).$$
These functions are defined on $(C_l\times C_l)\setminus diag(C_l\times C_l)$ for each fixed $u,v\in \mathbb{R}^n,$ where $diag(C_l\times C_l)=\{(x',x'):x'\in C_l\}$ and they are $l$-cylindrical remaining Rarita-Schwinger functions. Consequently,
 $Q_k^{C_l}Cot_{l,k}'(x',y',u,v)=0$ . Furthermore for each $l$ they represent a kernel for the
operator $Q_k^{C_l}$.

Now we will establish some integral formulas associated with the remaining Rarita-Schwinger operators on cylinders.

\begin{theorem}\cite{LR}(Stokes' Theorem for $Q_k$ operators) Let $\Omega'$ and $\Omega$ be domains in $\mathbb{R}^n$ and suppose the closure of $\Omega$ lies in $\Omega'$. Further suppose the closure of $\Omega$ is compact and the boundary of $\Omega,$ $\partial\Omega$, is piecewise smooth. Then for $f, g \in C^1(\Omega',$$\mathcal{M}_{k-1})$, we have
\begin{multline*}
\di\int_\Omega[(g(x,u)uQ_{k,r}), uf(x,u))_u+(g(x,u)u, Q_kuf(x,u))_u]dx^n\\
=\di\int_{\partial\Omega}\left(g(x,u)u, (I-P_k)d\sigma_xuf(x,u)\right)_u\\
=\di\int_{\partial\Omega}\left(g(x,u)ud\sigma_x(I-P_{k,r}), uf(x,u)\right)_u.
\end{multline*}
Where $Q_{k,r}$ is the right remaining Rarita-Schwinger operator.
\end{theorem}

\begin{theorem}\cite{LR}(Borel-Pompeiu Theorem for $Q_k$ operators)Let $\Omega'$ and $\Omega$ be as in the previous Theorem. Then for $f\in C^1(\Omega',$$\mathcal{M}_{k-1})$ and $y\in \Omega,$ we obtain
\begin{multline*}
uf(y,u)=\di\int_\Omega(H_k(x-y,u,v),Q_kvf(x,v))_vdx^n\\
-\di\int_{\partial\Omega}\left(H_k(x-y,u,v), (I-P_k)d\sigma_xvf(x,v)\right)_v.
\end{multline*}
 \end{theorem}

Applying Stokes' Theorem and Borel-Pompeiu Theorem for the $Q_k$ operator in $\mathbb{R}^n$, we may have:

\begin{theorem}\label{BPT0} \quad Let $V$ be a bounded domain in $\mathbb{R}^n$ and $\overline{V}$ be
the closure of $V$. For each $x\in \overline{V},$ the shifted lattice $x+\mathbb{Z}^l$
intersected with $V$ only contains the points $x.$ Suppose that the boundary of $V$, $\partial {V}$, is piecewise
smooth and $\overline{V}$ is compact. Further suppose $g(x,u): \overline{V}\times \mathbb{R}^n\to Cl_n$ is a
monogenic homogeneous polynomial of degree $k-1$ in $u$ and with respect to $x$ is $C^1$. Then for $1\leq l\leq n-1$ and each $y\in V,$
\begin{multline*}
vg(y,v)=\di\int_{\partial V}\di\int_V(Cot_{l,k}(x-y,u,v),Q_kuf(x,u))_u\\
-\left(Cot_{l,k}(x-y,u,v),(I-P_k)d\sigma_xug(x,u)\right)_udx^n.
\end{multline*}
\end{theorem}

Now using the projection map $\pi_l$, we obtain the Borel-Pompeiu Theorem for the $Q_k$ operators on the cylinder $C_l.$

\begin{theorem}\label{BPT}\quad (Borel-Pompeiu Theorem for $Q_k^{C_l}$)
Suppose $V'$ is a domain in $C_l$ with compact closure and smooth boundary. Suppose $g(x,u)$ is defined as in Theorem \ref{BPT0}.
Then for $1\leq l\leq n-1$ and each $y'\in V',$
\begin{multline*}
vg'(y',v)=\di\int_{\partial V'}\left(Cot_{l,k}'(x',y',u,v), (I-P_k)d\sigma'_{x'}ug'(x',u)\right)_u\\
\qquad \quad -\di\int_{V'}(Cot_{l,k}'(x',y',u,v),Q_k^{C_l}ug'(x',u))_ud\mu(x'),
\end{multline*}
where $x'=\pi_l(x), d\sigma'_{x'}=\partial_x\pi_ld\sigma_{x},$ $\partial_x\pi_l$ is the derivative of $\pi_l$ at $x$,
 $\mu$ is the projection of Lebesgue measure on $\mathbb{R}^n$ onto $C_l,$ $d\sigma_{x}= n(x)d\sigma(x)$, $n(x)$ is the unit exterior normal vector at $x = \pi_l^{-1}(x')$ and $d\sigma(x)$ is the surface measure
element.
\end{theorem}

\begin{theorem}\quad (Cauchy integral formula for $Q_k^{C_l}$)
Suppose the hypotheses as in Theorem \ref{BPT0} and $ug'(x',u)$ is annihilated by
the operator $Q_k^{C_l},$ then for $1\leq l\leq n-1$ and each $y'\in V',$
we have $$vg'(y',v)=\di\int_{\partial V'}\left(Cot_{l,k}'(x',y',u,v), (I-P_k)d\sigma'_{x'}ug'(x',u)\right)_u.$$
\end{theorem}

If the function given in Borel-Pompeiu Theorem has its support with respect to $x'$
in $V'\subset C_l,$ then we have the following:
\begin{theorem}
\quad $\di\int_{V'}-(Cot_{l,k}'(x',y',u,v),Q_k^{C_l}u\psi (x',u))_udx'^n=v\psi(y',v),$ \\for $\psi\in C^{\infty}(C_l\times \mathbb{R}^n).$
\end{theorem}

 We may introduce a Cauchy transform for the remaining Rarita-Schwinger operator $Q_k^{C_l}:$
\begin{definition} For a subdomain $V'$ of cylinder $C_l$ and a function $g': V'\times \mathbb{R}^n \longrightarrow Cl_n,$ which is monogenic in $u$ with degree $k-1$, the Cauchy
transform of $f'$ is formally defined to be
$$
(Tvg')(y',v)=-\int_{V'} \left(Cot_{l,k}'(x',y',u,v), ug'(x',u)\right)_udx'^n,\, y'\in V'.
$$
\end{definition}

Consequently, one may obtain:
\begin{theorem}
\quad $Q_k^{C_l}\di\int_{V'}(Cot_{l,k}'(x',y',u,v),u\psi (x',u))_udx'^n=v\psi(y',v),$\\ for $\psi\in C^{\infty}(C_l\times \mathbb{R}^n).$
\end{theorem}

Similarly, we can carry on the theory of the remaining Rarita-Schwinger operators to the setting of conformally spinor bundles over the cylinders.

Junxia Li \quad Email:  jxl004@uark.edu\\
John Ryan \quad Email: jryan@uark.edu\\
Carmen J. Vanegas \quad Email: cvanegas@usb.ve


\begin{thebibliography}{BSSV1}


\bibitem [BDS]{BDS} F. Brackx, R. Delanghe R. and F. Sommen, \textit{Clifford Analysis.} Pitman, London, 1982.

\bibitem [BSSV]{BSSV} J. Bure\v{s}, F. Sommen, V. Sou\v{c}ek and P. Van Lancker, \textit{Rarita-Schwinger Type Operators in Clifford Analysis.}
J. Funct. Anal. \textbf{185} (2001), No.2, 425--455.

\bibitem [BSSV1]{BSSV1} J. Bure\v{s}, F. Sommen, V. Sou\v{c}ek and P. Van Lancker, \textit{Symmetric Analogues of Rarita-Schwinger Equations.}
Annals of Global Analysis and Geometry \textbf{21} (2002), 215--240.

\bibitem [DLRV]{DLRV} C. Dunkl,J. Li, J. Ryan and P. Van Lancker, \textit{Some Rarita-Schwinger Operators.} Submitted (2011).

\bibitem [E]{E} G. Eisenstein, \textit{Genaue Un\-ter\-suchung der un\-endlichen Doppelproducte, aus welchen die elliptischen
Functionen als Quotienten zusammengesetzt sind, und der mit ihnen zusammenh\"{a}ngenden Doppelreihen
(als eine neue Begr\"{u}ndung der Theorie der elliptischen Functionen mit besonderer
Ber\"{u}cksichtigung ihrer Analogie zu den Kreisfunctionen).}
Crelle's Journal \textbf{35}, (1847), 153--274.

\bibitem [K]{K} R. S. Krausshar, \textit{Generalized Analytic Automorphic Forms in Hypercomplex Spaces.}
Frontiers in Mathematics, Birkh\"{a}user Verlag, 2004.

\bibitem[K1]{K1} R. S. Krausshar, \textit{Monogenic multiperiodic functions in Clifford analysis.} Complex Variables Theory Appl.
\textbf{46}, (2001), No.4, 337--368.

\bibitem [KR]{KR} R. S. Krausshar and J. Ryan, \textit{Clifford and Harmonic analysis on Cylinders and Tori.}
Rev. Mat. Iberoamericana \textbf{21}, (2005), No.1, 87--110.

\bibitem[KR1]{KR1} R. S. Krausshar and J. Ryan,\textit{Some conformally flat spin manifolds, Dirac operators and automorphic forms.}
 J. Math. Anal. Appl. \textbf{325}, (2007),  No.1, 359–-376.

\bibitem [LR]{LR} J. Li and J. Ryan, \textit{Some Operators Associated to Rarita-Schwinger Type Operators.} Complex Variables and Elliptic Equations, Accepted October (2011).

\bibitem [LRV]{LRV} J. Li, J. Ryan and C. J. Vanegas, \textit{Rarita-Schwinger Type Operators on Spheres and Real Projective Space.} Submitted (2011).

\bibitem [MP]{MP} R. Miatello, R. Podesta, \textit{Spin structures and spectra of $\mathbb{Z}_2$ manifolds.}
Math. Z. \textbf{247}, (2004), 319--335.

\bibitem[Va]{Va}  P. Van Lancker, \textit{Higher Spin Fields on Smooth Domains,} in Clifford Analysis and Its Applications, Eds. F. Brackx, J.S.R. Chisholm and V. Sou\v{c}ek, Kluwer, Dordrecht (2001), 389-398.

\bibitem[Va1]{Va1} P. Van Lancker, \textit{Rarita-Schwinger Fields in the Half Space.} Complex Variables and Elliptic Equations,
\textbf{51}, (2006), 563-579.

\end{thebibliography}
\end{document}